Cheng Lee
chengli@uab.edu

# A Bivariate Competing-Risks Model with One Termination Event

**Abstract**
The likelihood function for a competing-risks model with one fatal and one non-fatal event is proposed. A bivariate Weibull using the likelihood function is applied to the Stanford Heart Transplant Data.

**Keywords**
Competing-risks, Survival analysis, Bivariate Weibull

**Introduction**
Competing-risks models in survival analysis have been studied by several authors. For example, Kalbfleisch & Prentics (1980) and Lawless (1982) have discussed multivariate survival models with independent and dependent risks. Most of time, only the minimum of the survival time of the events can be observed. This is true when every event is the termination event that is every event is fatal and when any of the fatal events occurs to the individual in the study, the individual is uncensored due to the event and will no longer stay in the study. However, for a bivariate competing-risks model, when one event is non-fatal and the other is fatal, the individual who first experiences the non-fatal event still stays in the study until the occurrence of the fatal event or the end of the study. To accommodate this type of survival analysis, the likelihood function of the one-fatal-and-one-non-fatal competing-risks model is therefore proposed. A bivariate Weibull model with the proposed likelihood function is applied to the Stanford Heart Transplant Data.

**The likelihood function**
Let $X$ and $Y$ be the random variables of survival time or failure to, respectively, event $A$ and $B$. The likelihood function shown by Lawless (1982) for a bivariate competing-risks model with both $X$ and $Y$ observable is

$$L = \prod_{i=1}^{N} f_{XY}\left(t_{x_i}, t_{y_i}\right)^{p_i} \times \left(\left[-\frac{\partial}{\partial x} S_{XY}(x,y)\right]_{x=t_{x_i}, y=t_{y_i}}\right)^{q_i}$$

$$\times \left(\left[-\frac{\partial}{\partial y} S_{XY}(x,y)\right]_{x=t_{x_i}, y=t_{y_i}}\right)^{r_i} \times S_{XY}\left(t_{x_i}, t_{y_i}\right)^{1-p_i-q_i-r_i} \qquad (1)$$

where $f_{XY}(.)$ and $S_{XY}(.)$ denote, respectively, the joint probability density function and the joint survival function of $X$ and $Y$, and $p_i$, $q_i$, and $r_i$ are event indicators. In the likelihood function, the first factor is said to be uncensored due to event $A$ at $t_{x_i}$ and also uncensored due to event $B$ at $t_{y_i}$. The second factor is said to be uncensored due to event $A$ at $t_{x_i}$ and censored due to event $B$ at $t_{y_i}$. The third factor is said to be censored due to event $A$ at $t_{x_i}$ and uncensored due to event $B$ at $t_{y_i}$. The fourth factor is said to be censored due to event $A$ at $t_{x_i}$ and also censored due to event $B$ at $t_{y_i}$ either during or at the end of the study.

In the case when event $B$ is a termination event, then, event $A$ must always occur before event $B$ or never occurs when event $B$ occurs first because event $B$ is fatal. Suppose the study ends at time $t$ that is every individual is censored at time $t$, then there exists the probability of $X < Y$. However, the probability of $X > Y$ equals zero after the end of the study because of the fatal event $B$. The likelihood function for such a case becomes

$$\prod_{i=1}^{N} f_{XY}\left(t_{x_i}, t_{y_i}\right)^{p_i} \left(\left[-\frac{\partial}{\partial x} S_{XY}(x,y)\right]_{x=t_{x_i}, y=t_{y_i}}\right)^{q_i} f_Y\left(t_{y_i}\right)^{r_i}$$

$$\times \left(S_{XY}(t_i, t_i) + \int_{t_i}^{\infty} \left(\frac{\partial}{\partial y} S_{XY}(x,y)\right)_{x=y} dy\right)^{1-p_i-q_i-r_i} \quad .(2)$$

The first and the second factor are the same in Lawless' likelihood function. However, the third factor of likelihood function is the marginal probability density of $Y$ because the probability does not exist for the occurrence of event $A$ after event $B$ when an individual is uncensored due to event $B$ at $t_{y_i}$. The fourth factor is the probability of $t_i < X < Y$ where $t_i$ is the censoring time for both event $A$ and event $B$. The fourth factor is expressed by two terms. The first term is the joint survival function evaluated at $t_i$. The second term is the integration from $t_i$ to the infinity over $y$ on the integrand that is the derivative of the joint survival function with respect to $y$ then replaced $x$ by $y$. The derivation of the fourth factor is shown in the Appendix.

**Application**
The proposed likelihood function is applied to the Stanford Heart Transplant Survival Data in Crowley & Hu (1977) with no covariates. Noura & Read (1990) has analyzed the same data using the Weiubll model assuming the logarithm of the cumulative hazard function is linear to the logarithm of the survival time. With the same assumption, the bivariate Weibull model by Lu and Bhattacharyya (1990) is adopted in this application. The joint survival function of the bivariate Weibull mdoel of $X$ and $Y$ proposed by them is

$$S_{XY}(x,y) = Exp\left\{-\left[\left(\frac{x}{\lambda_1}\right)^{\frac{\gamma_1}{\alpha}} + \left(\frac{y}{\lambda_2}\right)^{\frac{\gamma_2}{\alpha}}\right]^{\alpha}\right\}$$ where $\lambda_1, \gamma_1, \lambda_2$ and $\gamma_2$ are positive and $0 < \alpha \leq$

0. The parameter $\alpha$ measures the association of $X$ and $Y$ and they are independent when $\alpha$ is equal to 1. In this application, $X$ is the time to heart transplant and $Y$ is the time to death either with or without the heart transplant.

Table 1 shows the 4 categories of survival time for $X$ and $Y$ corresponding to the 4 factors in the likelihood function. In the data, 43 individuals experienced both event $A$ and $B$; 24 individuals received transplants and were alive when censored; 29 individuals died before receiving transplants and 4 individuals did not receive transplants and were alive when censored during or at the end of the study. Records of individuals with $x$ equal to zero or $y$ equal to zero were deleted.

The estimates of parameters and their standard deviations based on the second derivatives valued at the maximized log-likelihood function are in Table 2. The shape parameters of $X$ and $Y$ are both negative indicating both have a decreasing hazard function. The mean denoted as $E(.)$, variance denoted as $Var(.)$, and the correlation denoted as $Corr(.)$ of $X$ and $Y$ using formulae by Lu and Bhattacharyya (1990) are in Table 3.

From Table 3, the correlation between $X$ and $Y$ is 0.3406 indicating that waiting time to heart transplant is not much correlated with the time to death.

## Appendix

The probability that $t < x < y < \infty$ is

$$\Pr(t < x < y < \infty)$$

$$= \int_t^\infty \int_t^y f_{XY}(x,y) \, dx \, dy$$

$$= \int_t^\infty \left[ \int_t^\infty f_{XY}(x,y) \, dx - \int_y^\infty f_{XY}(x,y) \, dx \right] dy$$

$$= S_{XY}(t,t) - \int_t^\infty \int_y^\infty f_{XY}(x,y) \, dx \, dy$$

Considering $\int_y^\infty f_{XY}(x,y) \, dx$,

$$\int_y^\infty f_{XY}(x,y) \, dx$$

$$= \int_0^\infty f_{XY}(x,y) \, dx - \int_0^y f_{XY}(x,y) \, dx$$

$$= f_Y(y) - \int_0^y \left( \frac{\partial}{\partial x} \left( \frac{\partial}{\partial y} F_{XY}(x,y) \right) \right) dx$$

$$= \frac{\partial}{\partial y} F_Y(y) - \left[ \frac{\partial}{\partial y} F_{XY}(x,y) \right]_{x=y} \quad (x=y \text{ denotes that } x \text{ is replaced by } y)$$

$$= \left[ \frac{\partial}{\partial y} (F_Y(y) - F_{XY}(x,y)) \right]_{x=y}$$

$$= \left[ \frac{\partial}{\partial y} (-S_{XY}(x,y)) \right]_{x=y}$$

Note that $S_{XY}(x,y) = 1 - F_X(x) - F_Y(y) + F_{XY}(x,y)$ (Li (1997))

Therefore,

$$\Pr(t < x < y < \infty)$$

$$= S_{XY}(t,t) + \int_t^\infty \left[ \frac{\partial}{\partial y} S_{XY}(x,y) \right]_{x=y} dy$$

## References

Crowley, J. & Hu, M. (1977). Covariance Analysis of Heart Transplant Survival Data. Journal of The American Statistical Association, Vol. 72, No. 357, pp 27-36

Kalbfleisch, J. D. & Prentice, R. L. (1980). The Statistical Analysis of Failure Time Data. John Wiley and Sons, New York, USA.
Lawless, J. F. (1982) Statistical Models and Methods for Lifetime Data, John Wiley and Sons, New York, USA.

Li, C. L. (1997). A Model for Informative Censoring, Ph.D. dissertation, The University of Alabama at Birmingham

Lu, J-C. & Bhattacharyya, G. (1990) Some New Constructions of Bivariate Weibull Models. *Annals of The Institute of Statistical Mathematics* 42, 543-559.

Noura, A. A. & Read, K. L. Q. (1990) Proportional Hazards Changepoint Models in Survival Analysis. Applied Statistics, Vol. 39, No. 2, pp 241-253.

Table 1

| p | q | r | x | Y | Transplant Status | Survival Status |
|---|---|---|---|---|---|---|
| 1 | 0 | 0 | Time to transplant | Time to death | Received | Dead |
| 0 | 1 | 0 | Time to transplant | Censoring time due to death | Received | Alive when censored |
| 0 | 0 | 1 | Time to death | Time to death | Not received | Dead |
| 0 | 0 | 0 | Time to censoring* | Time to censoring* | Not received | Alive when censored |

* Time to censoring is the time to the minimum of the date last seen and the end of the study, April 1, 1974

Table 2

| Parameter | Estimate | Standard Deviation |
|---|---|---|
| $\alpha$ | 0.5596 | 0.07356 |
| $\lambda_1$ | 35.5837 | 7.5053 |
| $\gamma_1$ | 0.5587 | 0.04529 |
| $\lambda_2$ | 385.6361 | 92.2737 |
| $\gamma_2$ | 0.48300 | 0.04849 |

Table 3

| Statistics | Estimate |
|---|---|
| $E(X)$ | 59.1578 |
| $E(Y)$ | 823.8204 |
| $Var(X)$ | 12972.4236 |
| $Var(Y)$ | 3742873.8509 |
| $Corr(X, Y)$ | 0.3406 |